\documentclass{amsart}
\usepackage{amsmath}
\usepackage{amsfonts}
\usepackage{amssymb}

\theoremstyle{theorem}

\theoremstyle{theorem}
\newtheorem{corollary}{Corollary}
\theoremstyle{theorem}
\newtheorem{lemma}{Lemma}
\theoremstyle{definition}
\newtheorem{definition}{Definition}
\theoremstyle{remark}
\newtheorem{remark}{Remark}
\theoremstyle{theorem}
\newtheorem{theorem}{Theorem}
\theoremstyle{remark}
\newtheorem{example}{Example}
\theoremstyle{theorem}
\newtheorem{algorithm}{Algorithm}

\DeclareMathOperator{\rank}{rank}
\DeclareMathOperator{\ord}{ord}
\DeclareMathOperator{\charact}{char}
\DeclareMathOperator{\Quot}{Quot}

\begin{document}

\title[Computation of Characteristic Sets]{On Computation of Kolchin Characteristic Sets: Ordinary and Partial Cases}
\author{Marina Kondratieva and Alexey Ovchinnikov}
\email{kondra\_m@shade.msu.ru, aiovchin@ncsu.edu}
\address{Moscow State University, Department of Mechanics and Mathematics\\
North Carolina State University, Department of Mathematics}
\date{December 4, 2004}
\thanks{The work was partially supported by the Russian Foundation for Basic Research, project no. 02-01-01033
and by NSF Grant CCR-0096842.}

\begin{abstract}
In this paper we study the problem of computing a Kolchin characteristic set of a radical differential ideal. The central part of the article is the presentation of algorithms solving this problem in two principal cases: for ordinary differential polynomials and in the partial differential case. Our computations are mainly performed with respect to orderly rankings. We also discuss the usefulness of regular and characteristic decompositions of radical differential ideals. In the partial differential case we give an algorithm for computing characteristic sets in the {\it special} case of radical differential ideals satisfying the property of consistency. For this class of ideals we show how to deal with arbitrary differential rankings.
\end{abstract}

\maketitle

\newcommand{\Le}{\leqslant}
\newcommand{\Ge}{\geqslant}

\begin{section}{Introduction}
This paper is devoted to study radical differential ideals and their characteristic sets. The concept of a characteristic set introduced by Ritt and Kolchin is one of the most important notions in differential algebra. The problem of computing a characteristic set of a radical differential ideal represented by a finite set of its generators is not completely solved yet especially in the partial differential case. In the case of ideals
in rings of polynomials in a finite number of variables this problem was studied and completely solved by Gallo and Mishra \cite{Gal1,Gal2,Gal3}. It was also investigated by Aubry, Lazard and Moreno Maza \cite{Tri}.

So, it is very natural to study this problem in rings of differential polynomials.
The most important contributions of this article are:
\begin{itemize}
\item an algorithm for computing characteristic sets of arbitrary radical differential ideals w.r.t. orderly rankings in the {\it ordinary} case;
\item an algorithm for computing characteristic sets of radical differential ideals satisfying the property of consistency in the {\it partial} differential case. We do this w.r.t. both orderly (Sections~\ref{Ord1} and \ref{Ord2}) and arbitrary  (Section~\ref{Arbitrary}) rankings.
\end{itemize}
We also conjecture a method of solving this problem in non-ordinary cases for arbitrary radical differential ideals (see Section~\ref{Conjecture}).

We use other techniques and methods than those used by Gallo and Mishra. However, their algorithm for computing a characteristic set of an algebraic ideal plays an important role in Algorithm~\ref{CharOrd} and Algorithm~\ref{CharAlg} of this paper. The methods developed by Sadik \cite{Sadik} help us to obtain several bounds for characteristic sets w.r.t. different rankings. 

In \cite{Sadik2} Sadik constructed an algorithm for computing a Kolchin
characteristic set of a radical differential ideal w.r.t. an {\it elimination} ranking in the {\it ordinary} case. One can also get some bounds for orderly rankings from \cite{Sadik2}, but these bounds are bigger then ours (see Remark~\ref{ComparisonSadik}). 
Thus, the results of this paper (Theorem~\ref{OrdinaryBound} and Theorem~\ref{Compute}) are {\it new} in comparison with Sadik's ones, also because Sadik did not give any algorithm for the partial differential case but we do.

Ten years ago a technique for effective and factorization-free computations in the radical differential ideal theory was developed by Boulier, Lazard, Ollivier and Petitot (see \cite{Bou1} and \cite{Bou2}). In \cite{Fac, Dif} Hubert continued to develop this problem and introduced the notions of {\it characterizable} ideal and {\it characteristic decomposition} of a radical differential ideal. This decomposition of the ideal helps us to solve many problems concerning the system of differential equations associated with the ideal and to test the radical membership.

It should be emphasized that a characteristic decomposition of a radical differential ideal does not give us full information about the ideal.
In some important cases a representation of this ideal by characteristic components cannot replace a representation of the ideal
by its generators as a radical differential ideal.
For example, at this moment one cannot check the inclusion of a radical differential ideal to another radical
differential ideal knowing only a characteristic decomposition of the first one (see \cite{Kon,Ov1,Ov2}). This problem is closely related to the well-known Ritt problem.
In this case it is necessary to know generators of the ideal and characteristic decomposition is partially useless.

Nevertheless, in this paper we show that a characteristic decomposition of a radical differential ideal yields a lot of information about the ideal. Indeed,
characteristic decomposition allows us not only to test membership to this ideal but also to compute its characteristic set in Kolchin's sense.
Hence, the main contribution of this paper can be also considered as another application of a characteristic decomposition. Such a decomposition tells us a lot about a system of partial differential equations and we study what else one can do using a characteristic decomposition.

We also study the case of {\it partial} differential polynomials and give an algorithm (Algorithm~\ref{CharAlg}) for computing a Kolchin characteristic set. Theorem~\ref{Compute} in Section~\ref{Result} provides a theoretical basis for this. Radical differential ideals, for which we propose the algorithm, satisfy the property of consistency (see Definition~\ref{Consistency}) w.r.t. an orderly ranking. We will see that these ideals are {\it better} than arbitrary radical differential ideals in the {\it computational sense}.

In summary, although a characteristic decomposition cannot replace the ideal in computational sense, this decomposition allows us to compute such an important subset of a radical differential ideal as its characteristic set.
\end{section}

\begin{section}{Preliminaries}

\begin{subsection}{Basic definitions}
Differential algebra deals with differential rings and fields. These are commutative rings with $1$ and a basic set of differentiations $\Delta = \{\delta_1,\ldots,\delta_n\}$ on the ring. The case of $\Delta = \{\delta\}$ is called {\it ordinary}. If $R$ is an ordinary differential ring and $y \in R$ we denote $\delta^ky$ by $y^{(k)}$. The ring of differential polynomials was introduced to deal with algebraic differential equations.

Recent tutorials on constructive differential ideal theory are presented in \cite{Dif, Sit}. We also use the Gr\"obner bases technique discussed in detail in \cite{Wei}. The definition of the ring of {\it differential polynomials} in $l$ variables over a differential field $k$ is given in \cite{Kol, Pan, Rit}. This ring is denoted by $k\{y_1,\dots,y_l\}.$ We consider the case of $\charact k=0$ only. We denote polynomials by $f, g, h, \ldots$ and use the notation $I, J, P, Q$ for ideals.

We need the notion of reduction for algorithmic computations. First, we introduce a {\it ranking} on the set of differential variables of $k\{y_1,\ldots,y_l\}$. Construct the multiplicative monoid $\Theta = (\delta_1^{k_1}\delta_2^{k_2}\cdots\delta_n^{k_n}, k_i \Ge 0).$ The ranking is a total ordering on the set $\{\theta y_i\}$ for each $\theta \in \Theta$ and $1\Le i \Le l$
satisfying the following conditions:
\begin{enumerate}
\item $\theta u \Ge u,$
\item $u \Ge v \Longrightarrow \theta u \Ge \theta v.$
\end{enumerate}
In later discussions we suppose that a ranking is fixed.

Let $u$ be a differential variable in $k\{y_1,\ldots,y_l\}$, that is, $u = \theta y_j$ for a differential operator $\theta = \delta_1^{k_1}\delta_2^{k_2}\cdots\delta_n^{k_n} \in \Theta$ and $1\Le j \Le l$. A ranking is said to be {\it orderly} iff $\ord u > \ord v$ implies $u > v$ for all differential variables $u$ and $v$.
A ranking $>_{el}$ is called {\it elimination} iff $y_i >_{el} y_j$ implies $\theta_1y_i >_{el} \theta_2y_j$ for all $\theta_1, \theta_2 \in \Theta$.

The highest ranked derivative $\theta y_j$ appearing in a differential polynomial $f \in k\{y_1,\dots,y_l\} \setminus k$ is called the leader of $f$. We denote the leader by $u_f$. Represent $f$ as a univariate polynomial in $u_f$:
$$
f = I_f u_f^n + a_1 u_f^{n-1} + \ldots + a_n.
$$
The polynomial $I_f$ is called the {\it initial} of $f$.

Apply any $\delta \in \Delta$ to $f$:
$$ \delta f = \frac{\partial f}{\partial u_f}\delta u_f + \delta I_f u_f^n +
\delta a_1 u_f^{n-1}+\ldots + \delta a_n.
$$
The leading variable of $\delta f$ is $\delta u_f$ and the initial of $\delta f$ is called the {\it separant} of $f$. We denote it by $S_f$. Note that for all $\theta \in \Theta,$ $\theta \ne 1,$ each $\theta f$ has the initial equal to $S_f$.

Define the ranking on differential polynomials. We say that $f > g$ iff $u_f > u_g$ or in the case of $u_f = u_g$
we have $\deg_{u_g} f > \deg_{u_g} g$. Let $F \subset k\{y_1,\dots,y_l\}$ be a set of differential polynomials. For the differential and radical differential ideal generated by $F$ in $k\{y_1,\dots,y_l\}$, we use the notation $[F]$ and $\{F\}$, respectively.

We say that a differential polynomial $f$ is {\it partially reduced} w.r.t. $g$ iff no proper derivative of $u_g$ appears in $f$. A differential polynomial $f$ is {\it reduced} w.r.t. $g$ iff $f$ is partially reduced w.r.t. $g$ and $\deg_{u_g} f < \deg_{u_g} g$. Consider any subset $\mathbb{A} \subset k\{y_1,\ldots,y_l\}$. We say that $\mathbb{A}$ is autoreduced iff $\mathbb{A} \cap k = \varnothing$ and each element of $\mathbb{A}$ is reduced w.r.t. all the others. Every autoreduced set is finite (see \cite[Chapter I, Section 9]{Kol}). For autoreduced sets we use capital letters $\mathbb{A, B, C,}$ \ldots.

We denote the product of the initials and the separants of the elements of $\mathbb{A}$ by $I_\mathbb{A}$ and $S_\mathbb{A}$, respectively. Denote $I_\mathbb{A}\cdot S_\mathbb{A}$ by $H_\mathbb{A}$. Let $S$ be a finite set of differential polynomials. Denote by $S^\infty$ the multiplicative set containing $1$ and generated by $S$. Let $I$ be an ideal in a commutative ring $R$. Let $I:S^\infty = \{a \in R\:|\:\exists s \in S^\infty: sa \in I\}$.
If $I$ is a differential ideal then $I:S^\infty$ is also a differential ideal (see \cite{ Kol, Rit, Pan, Sit}).

If we want to enumerate the elements of $\mathbb{A}$ we write the following: $\mathbb{A} = A_1, A_2, \ldots, A_p$.
Let $\mathbb{A} = A_1,\ldots,A_r$ and $\mathbb{B} = B_1,\ldots,B_s$ be autoreduced sets. Let the elements of $\mathbb{A}$ and $\mathbb{B}$ be arranged in order of increasing rank. We say that $\mathbb{A}$ has lower rank than $\mathbb{B}$ iff there exists $k \Le r, s$ such that $\rank A_i$ = $\rank B_i$ for $1 \Le i < k$ and $\rank A_k < \rank B_k$, or if $r > s$ and $\rank A_i = \rank B_i$ for $1 \Le i \Le s$. We say that $\rank\mathbb{A} = \rank\mathbb{B}$ iff $r=s$ and $\rank A_i = \rank B_i$ for $1 \Le i \Le r$.

Consider two differential polynomials $f$ and $g$ in $R = k\{y_1,\ldots, y_l\}$. Let $I$ be the differential ideal in $R$ generated by $g$. Applying a finite number of differentiations and pseudo-divisions one can compute a {\it differential partial remainder} $f_1$ and a {\it differential remainder} $f_2$
of $f$ w.r.t. $g$ such that there exist $s \in S_g^\infty$ and $h \in H_g^\infty$ satisfying $sf \equiv f_1$ and $hf \equiv f_2 \mod I$ with $f_1$ and $f_2$ partially reduced and reduced w.r.t. $g$, respectively (see \cite{Fac} for definitions and an algorithm for computing remainders).

Let $\mathbb{A}$ be an autoreduced set in $k\{y_1,\ldots,y_l\}$. Consider the polynomial ring $k[x_1,\ldots,x_n]$ with $x_1,\ldots,x_n$ belonging to $\Theta Y$ for $Y = y_1,\ldots,y_l$. Let $U, V \subset \{x_1,\ldots,x_n\}$ be the sets of ``leaders'' and ``non-leaders'' appearing in the autoreduced set $\mathbb{A}$, respectively. We denote $k[x_1,\ldots, x_n]$ by $k[V][U]$ and the leader of $A_i$ by $u_{A_i}$ or $u_i$ for each $1 \Le i \Le p$.

The notion of a characteristic set in {\it Kolchin's sense} in characteristic zero is crucial in our further discussions. This was first used by Ritt for prime differential ideals, but Kolchin introduced characteristic sets for arbitrary differential ideals.

\begin{definition} \cite[page 82]{Kol} An autoreduced set of the lowest rank in an ideal $I$ is
called a {\it characteristic set} of $I$.
\end{definition}

We call these sets {\it Kolchin characteristic sets} in order to omit confusions with different notions, e.g., in \cite{Fac,Dif} characteristic sets are used in Kolchin's sense and in some other senses.

As it is mentioned in \cite[Lemma 8, page 82]{Kol}, in characteristic zero $\mathbb{A}$ is a characteristic set of a proper differential ideal $I$ iff each element of $I$ reduces to zero w.r.t. $\mathbb{A}$. Consider the definition of a characterizable radical differential ideal.

\begin{definition}\cite[Definition 2.6]{Fac} A radical differential ideal $I$ in $k\{y_1,\ldots,y_l\}$ is
said to be {\it characterizable} iff there exists a characteristic set $\mathbb{A}$ of $I$ in Kolchin's sense
such that $I = [\mathbb{A}]:H_\mathbb{A}^\infty.$
\end{definition}

The following definition makes a bridge between differential and commutative algebra. Let $v$ be a derivative in $k\{y_1,\ldots,y_l\}.$ $\mathbb{A}_v$ is the set of the elements of $\mathbb{A}$ and their derivatives that have a leader ranking strictly lower than $v$.
\begin{definition}\cite[III.8]{Kol}
$\mathbb{A}$ is {\it coherent} iff whenever $A, B \in \mathbb{A}$ are such that $u_{A}$ and $u_{B}$ have a common derivative: $v = \psi u_{A} = \phi u_{B}$, then $S_{B}\psi A - S_{A}\phi B \in (\mathbb{A}_v):H_\mathbb{A}^\infty.$
\end{definition}

We emphasize that a characteristic set of a differential ideal is a coherent autoreduced set (see \cite{Kol,Rit, Pan, Sit}).
\end{subsection}

\begin{subsection}{Important assertions}
Consider several important results concerning radical differential ideals in rings of differential polynomials. The technique described in \cite{Fac, Kol} helps us to cover some properties of these ideals.

\begin{theorem}\cite[III.8, Lemma 5]{Kol}\label{Rosenfeld} Let $\mathbb{A}$ be a coherent autoreduced set in $k\{y_1,\ldots,y_l\}$. Suppose that a differential polynomial $g$ is partially reduced w.r.t. $\mathbb{A}$. Then $g \in [\mathbb{A}]:H_\mathbb{A}^\infty$ iff $g \in (\mathbb{A}):H_\mathbb{A}^\infty$.
\end{theorem}

Note that Theorem~\ref{Rosenfeld} is also known as {\it Rosenfeld's lemma}.

\begin{theorem}\cite[Theorem 3.2]{Fac}\label{T20} Let $\mathbb{A}$ be an autoreduced set of $k[V][U]$. If $1 \notin (\mathbb{A}):S_\mathbb{A}^\infty$ then any minimal prime of $(\mathbb{A}):S_\mathbb{A}^\infty$ admits the set of non-leaders of
$\mathbb{A}$, $V$, as a transcendence basis. More specially, any characteristic set of a minimal prime of $(\mathbb{A}):S_\mathbb{A}^\infty$ has the same set of leaders as $\mathbb{A}$.
\end{theorem}

\begin{theorem}\cite[Theorem 4.5]{Fac}\label{Prime} Let $\mathbb{A}$ be a coherent autoreduced set of $R = k\{y_1,\ldots,y_l\}$ such that $1 \notin [\mathbb{A}]:H_\mathbb{A}^\infty$. There is a one-to-one correspondence between the minimal primes of $(\mathbb{A}):H_\mathbb{A}^\infty$ in $k[V][U]$ and the essential prime components of $[\mathbb{A}]:H_\mathbb{A}^\infty$ in
$R$. Assume $\mathbb{C}_i$ is a characteristic set of a minimal prime of $(\mathbb{A}):H_\mathbb{A}^\infty$. Then $\mathbb{C}_i$ is the characteristic set of a single essential prime component of $[\mathbb{A}]:H_\mathbb{A}^\infty$ (and vice versa).
\end{theorem}

\begin{lemma}\label{LoseVar} Let $\mathbb{A} = A_1, \ldots, A_p$ be an autoreduced set in the ring $k[x_1,\ldots,x_m] = R$ and a characteristic set of $(\mathbb{A}):I_\mathbb{A}^\infty$. Suppose that a polynomial $f = a_mx_t^m + \ldots + a_0 \in R$ is reducible to zero w.r.t. $\mathbb{A}$ and the indeterminate $x_t$ does not appear in $A_i$ for each $1\Le i \Le p$. Then $a_j$ is reducible to zero w.r.t. $\mathbb{A}$ for all $0 \Le j \Le m$.
\end{lemma}
\begin{proof}
Since $f$ is reducible to zero w.r.t. $\mathbb{A}$, there exists $I \in  I_\mathbb{A}^\infty$ such that
$$
I\cdot f = \sum\limits_{i=1}^p g_iA_i.
$$

Let $g_i = \sum\limits_{j=1}^{t_i}h_j x_t^j$ for each $1 \Le i \Le p$. Thus, we have $I\cdot \sum\limits_{k=0}^m a_k x_t^k = \sum\limits_{k=0}^q d_k x_t^k$ with $d_k \in (A_1,\ldots,A_p)$. Hence, $I\cdot a_i \in (\mathbb{A})$ for each $1 \Le i \Le m$,
that is, $a_i \in (\mathbb{A}):I_\mathbb{A}^\infty$. Since  $\mathbb{A}$ is a characteristic set of $(\mathbb{A}):I_\mathbb{A}^\infty$, we have that all $a_i$ are reducible to zero w.r.t. $\mathbb{A}$.
\end{proof}

\end{subsection}

\end{section}

\begin{section}{The ordinary case}\label{Result}
\begin{subsection}{Bounds for the orders of characteristic sets}
Denote \cite[Algorithm 7.1]{Fac} by {\sf $\chi$-Decomposition}. An input of this algorithm is a finite set of differential polynomials $F$ and
its output is a set $\mathfrak{C} = \mathbb{C}_1,\ldots,\mathbb{C}_n$ of characteristic sets $\mathbb{C}_i$ of characterizable ideals $[\mathbb{C}_i]:H_{\mathbb{C}_i}^\infty$ forming a characteristic decomposition of the radical differential ideal $\{F\}$ in $k\{y_1,\ldots,y_l\}$:
$$
\{F\} = [\mathbb{C}_1]:H_{\mathbb{C}_1}^\infty\cap\ldots\cap[\mathbb{C}_n]:H_{\mathbb{C}_n}^\infty.
$$
We will show how this decomposition helps to compute a characteristic set of $\{F\}$.

The main idea of Algorithm~\ref{CharOrd} is to move our problem into commutative algebra.
In order to do this we need to know bounds for the orders of elements of a characteristic set of a
radical differential ideal.

Let $R = k\{y_1,\ldots,y_l\}$ with $\Delta = \{\delta\}$. So, we are in the ordinary case. {\it Differential dimension} of a prime differential ideal $P$ is the maximal number $q$ such that $P \cap k\{y_{i_1},\ldots,y_{i_q}\} = \{0\}$.
If $f$ is a differential polynomial then $\ord f$ denotes the maximal order of differential variables appearing effectively in $f$. Let $\mathbb{A} = A_1,\ldots,A_p$ be an autoreduced set. Define the order of $\mathbb{A}$ by the following equality: $\ord\mathbb{A} = \ord A_1 + \ldots + \ord A_p$. If a set $\mathbb{C}$ is characteristic of the ideal $P$ w.r.t. an orderly ranking then by definition the order of the ideal $P$ equals $\ord\mathbb{C}$.

Denote by $P(s)$ the elements of $P$ of the order less than of equal to $s$. The set $P(s)$ is a prime ideal in the correspondent polynomial ring. As it is proved in \cite{Pan} the dimension of $P(s)$ is a polynomial in $s$ for $s \Ge h = \ord P$. More precisely, $\dim P(s) = q(s+1) + h$, where $q$ is the differential dimension of the ideal $P$. Moreover, $q = l-p$ and $p$ is the number of elements of a characteristic set of the ideal $P$
w.r.t. an orderly ranking. Thus, the number $p$ does not depend on an orderly ranking.

\begin{remark}
Almost all recent results in the dimensional theory for prime differential ideals (the theory of differential dimension polynomials) are presented in \cite{Pan}.
\end{remark}

\begin{lemma}\cite[Proposition 17]{Sadik}\label{OrderLess} Consider a prime differential ideal $P$ of differential dimension $q$ and of order $h$. For every subset $\{y_{i_1},\ldots,y_{i_{q+1}}\}$ of $\{y_1,\ldots,y_l\}$ the ideal $P$ contains a differential polynomial in the indeterminates $\{y_{i_1},\ldots,y_{i_{q+1}}\}$ with order less than or equal to $h$.
\end{lemma}

\begin{lemma}\cite[Lemma 23]{Sadik}\label{Separants} Let $P$ be a prime  differential ideal, then $P$ admits a characteristic set $\mathbb{C} = C_1,\ldots,C_p$ such that $\dfrac{\partial C_i}{\partial v_k^{({l_{i,k}})}}$ does not lie in $P$, where $v_k$ is a non-leader
and $l_{i,k} = \ord(\mathbb{C}_i, v_k)$ for a pair $(i,k)$ in $\{1,\ldots,p\}\times\{1,\ldots,l-p\}$.
\end{lemma}

A characteristic set of a prime differential ideal is not unique, e.g., consider the ideal $[x] \in k\{x,y\}$ and the elimination ranking with $x > y$. Then the set $y^{(n)}x$ is a characteristic set of the ideal $[x]$ for any $n \Ge 0$. Hence, if we do not propose any restriction, we will not be able to guarantee that a characteristic set of a prime differential ideal has restrictions for the orders of its characteristic sets.

In order to omit this problem, in \cite{Sadik} Sadik used {\it irreducible} characteristic set and proved in \cite[Lemma 19]{Sadik} that any prime differential ideal has an irreducible characteristic set. 

The following result will give us an opportunity to find a characteristic set of a prime differential ideal with good restrictions for the orders of its elements w.r.t. {\it any} differential ranking.

\begin{lemma}\label{IrreducibleSeparants} A prime differential ideal $P$ in $k\{y_1,\ldots,y_l\}$ admits a characteristic set $\mathbb{C} = C_1,\ldots,C_p$ with the following properties:
\begin{enumerate}
\item It is irreducible, that is, 
{\begin{enumerate}
\item $C_1$ is an irreducible polynomial in $k\{y_1,\ldots,y_l\}$,
\item each $C_{i+1}$ is irreducible in the ring $\Quot(k[V]/(C_1,\ldots,C_i):I_i^\infty)[U]$, where 
{\begin{itemize}
\item $V$ is the set of all variables appearing in the polynomials $C_1,\ldots,C_i$, 
\item $I_i^\infty$ is the multiplicative system generated by the initials of the polynomials $C_1,\ldots,C_i$,
\item $U$ is the set of that variables from $C_{i+1}$ that are not in $V$.
\end{itemize}}
\end{enumerate}
}
\item Let $y_t^{(s)}$ be a differential variable of the maximal order $s$ appearing in $\mathbb{C}$. Let also $y_t^{(s)}$ does not appear in $C_1,\ldots,C_{i-1}$ but does appear in $C_i$. Then $S_{i,t} = \dfrac{\partial C_i}{\partial y_t^{(s)}} \notin P$.
\end{enumerate}  
\end{lemma}
\begin{proof} This is just a combination of the proofs of Lemma~\ref{Separants} and \cite[Lemma 19]{Sadik}. Indeed,
let us take an irreducible characteristic set $\mathbb{C}$ given
by \cite[Lemma 19]{Sadik}. Let $y_t^{(s)}$ be a differential variable of the maximal order $s$ appearing in $\mathbb{C}$. Let also $y_t^{(s)}$ does not appear in $C_1,\ldots,C_{i-1}$ but does appear in $C_i$. 

If $S_{i,t} \in P$. Then replace $C_i$ by $S_{i,t}$ in $\mathbb{C}$ and apply the process producing irreducibility from \cite[Lemma 19]{Sadik} to the new $\mathbb{C}$. The polynomials  $C_{i+1},\ldots,C_p$ will not be affected, since we do not change the ideal over which we factor. Thus, we can proceed by induction on the degree of $y_t^{(s)}$ in $C_i$ and then by the maximal number $s$.
\end{proof}

\begin{theorem}\label{OrderlyElimination} Let $P$ be a prime differential ideal of order $h$ in $k\{y_1,\ldots,y_l\}$ and $>$ be a differential ranking. Then there exists a characteristic set $\mathbb{C} = C_1,\ldots,C_n$ of the ideal $P$ w.r.t the ranking $>$ such that the order in $y_t$ of each $C_i$ does not exceed $h$ for all $1\Le t \Le l$.
\end{theorem}
\begin{proof}
The proof of \cite[Theorem 24]{Sadik} will be valid for our purpose if we {\it modify} it in the following way. Denote the set $\{y_k\:|\: y_k$ is not a leader of any $C_j,$ $1 \Le j \Le n\}$ by $\mathfrak{N}.$
If for some $\theta \in \Theta$ and $t,$ $1 \Le t \Le l$ the variable $\theta y_t$ is the leader of some $C_j$ then $\ord(C_q,y_t) \Le h$
for all $1\Le q \Le n$ immediately by Lemma~\ref{OrderLess}, since $\mathbb{C}$ is autoreduced  and $\ord\theta \Le h$. If not, we just take a polynomial $f \in k\{y_{C_j},\mathfrak{N}\}\cap P$ of
order not greater than $h$. We have $f$ is reduced w.r.t. $\mathbb{C}$. Contradiction. 

Let $y_t \in \mathfrak{N}$ and $\mathbb{C}$ be a characteristic set given by Lemma~\ref{IrreducibleSeparants}. 
Now, the main idea is to reduce the polynomial $f_j \in k\{y_{C_j},\mathfrak{N}\}$ given by Lemma~\ref{OrderLess} w.r.t. $\mathbb{C}$ where $\theta y_{C_j}$ is the leader of $C_j$ for some $\theta \in \Theta$. Suppose that for some $j, 1\Le j \Le n,$ we have $\ord(C_j, y_t) > h$. 

Note that $\ord(f_j,y_{C_j}) \Ge \ord(C_j,y_{C_j})$. Indeed, we may choose $f_j$ with $I_{f_j} \notin P$. Hence, if $\ord(f_j,y_{C_j}) < \ord(C_j,y_{C_j})$ then the ideal $P$ would contain an element that is not reducible w.r.t. the characteristic set $\mathbb{C}$, because the ideal $P$ is prime. Contradiction.
Let $\ord(f_j,y_{C_j}) \Ge \ord(C_j,y_{C_j})$.

Our {\it improvement} is to choose $$\tilde{\mathbb{C}} = \arg\max\limits_{C_j \in \mathbb{C}}\ord(C_j,y_t)$$ and then take $C_i \in \tilde{\mathbb{C}}$ of the {\it lowest} possible rank instead of using induction as Sadik did.
Let $u_i = \theta_i y_i$ for some $\theta_i \in \Theta$ and $u_i$ be the leader of $C_i$.
We have $s = \ord(C_i, y_t) > h$ and $r_f = \ord(f_i,y_i) \Ge \ord(C_i,y_i) = r_C$, where $f_i = f_i(y_i,\mathfrak{N}) = I_{f_i}\left(y_i^{(r_f)}\right)^{n_f} + a_1\left(y_i^{(r_f)}\right)^{n_f-1}+\ldots+a_{n_f}$. 

Let us reduce each term (coefficients $a_j$, initial $I_{f_i}$ and its leader $y_i^{(r_f)}$) of $f_i$ first by $C_i$. We need to differentiate $C_i$ $q$ times and get the remainder $\tilde{f}$ where $0  \Le q \Le r_f - r_C$. Remember that $f_i$ depends only on $y_i, \mathfrak{N},$ and their derivatives. Hence, applying further steps of reduction to the terms of $\tilde{f}$ w.r.t. other $C_j$ we need to differentiate them {\it less} than $q$ times if $C_j \in \tilde{\mathbb{C}}$ or {\it not greater} than $q$ times if $C_j \notin \tilde{\mathbb{C}}$. Indeed, the set $\mathbb{C}$ is autoreduced and the variables to reduce can come just from derivatives of variables from $C_i$. Moreover, $C_i$ has the smallest rank in $\tilde{\mathbb{C}}$.

That is why after we reduce all leaders of $\mathbb{C}$ from $f$ we get the polynomial depending effectively on $y_t^{(s+q)}$ and $s+q \Ge s$. 
Its initial w.r.t. $y_t^{(s+q)}$ is equal to $$I_{C_1}^{i_1}\cdot\ldots\cdot I_{C_n}^{i_n}\cdot S_{C_1}^{j_1}\cdot\ldots\cdot S_{C_n}^{j_n}\cdot\tilde{I_{f_i}}\cdot\left(\frac{\partial C_i}{\partial y_t^{(s)}}\right)^{n_f},$$
where $i_1,\ldots,i_n,j_1,\ldots,j_n \in \mathbb{Z}_{\Ge 0}$ and $\tilde{I_{f_i}}$ is the remainder of $I_{f_i}$ w.r.t. $\mathbb{C}$ in the case of $r_f > r_C.$ Moreover, in the case of $r_f = r_C$ we are in \cite[Lemma 20]{Sadik} because of our choice of $C_i$ and immediately get the inequality $\ord(f_i,y_t) \Ge \ord(C_i,y_t)$. Continue with the ``differential'' case. 

Remember that $P$ is a prime ideal. Hence, $I_{C_1}^{i_1}\cdot\ldots\cdot I_{C_n}^{i_n}\cdot S_{C_1}^{j_1}\cdot\ldots\cdot S_{C_n}^{j_n} \notin P$, because $I_{C_j}$ and $S_{C_j} \notin P$ for all $j,$  $1 \Le j \Le n.$
Moreover, $P = [\mathbb{C}]:H_\mathbb{C}^\infty$ and $\mathbb{C}$ is a characteristic set of $[\mathbb{C}]:H_\mathbb{C}^\infty$. Thus, according to Lemma~\ref{LoseVar} the polynomial $\frac{\partial C_i}{\partial y_t^{(s)}}$  is reducible to zero w.r.t. $\mathbb{C}$ that contradicts either to Lemma~\ref{IrreducibleSeparants}, or the discussions in the previous paragraph, because for a prime differential ideal the fact that an element is reducible to zero w.r.t. a characteristic set is nothing else the element belongs to the ideal.
\end{proof}

\begin{remark} Note that Theorem~\ref{OrderlyElimination} is actually a generalization of Sadik's result \cite[Theorem 24]{Sadik} that was proved just for elimination rankings.
\end{remark}

So, now we are ready to prove a final bound for a characteristic set of a radical differential ideal and to obtain an algorithm computing this set.

\begin{theorem}\label{OrdinaryBound}  Let $I$ be a radical differential ideal in $k\{y_1,\ldots,y_l\}$.
Let $I = \bigcap_{i=1}^n[\mathbb{C}_i]:H_{\mathbb{C}_i}^\infty$ be a characteristic decomposition w.r.t. an orderly ranking with $\mathbb{C}_i = C_i^1,\ldots,C_i^{p_i}$. Let $h$ be the maximal order of $\mathbb{C}_i$ for $1\Le i \Le n$. Then the lowest differentially autoreduced subset of an algebraic characteristic set of the ideal 
$$I' = \bigcap_{i=1}^n(\theta_i^jC_i^j, \ord\theta_i^j u_{C_i^j} \Le h):H_{\mathbb{C}_i}^\infty$$
is a characteristic set of $I$ w.r.t. the orderly ranking.
\end{theorem}
\begin{proof}
Let $I = \bigcap_{j=1}^a P_j$ be a minimal prime decomposition  and $I$ have a  characteristic set $\mathbb{C} = C_1,\ldots,C_p$ w.r.t. the orderly ranking. Suppose that for some $i,$ $1 \Le i \Le p$, we have $\ord (C_i, y_{j_i}) > h$, where $\theta y_{j_i} = u_i$ is the leader of $C_i$ and $\theta \in \Theta$. We may suppose this because the ranking is orderly. Note that $i \Ge 2$, because the leader of the differential polynomial  $C_1$ appears among the leaders of characteristic sets of the ideals $P_j$. If not then multiplying all the lowest elements of characteristic set of the ideals $P_j$ we get a non-reducible w.r.t. $\mathbb{C}$ element of $I$, that is a contradiction.

Let $I_i$ be the initial of $C_i$ with $I_i \in P_1,\ldots,P_t$ and $I_i \notin P_{t+1},\ldots,P_a$. Concentrate our attention at $P_{t+1},\ldots,P_a$. Denote a characteristic set of $P_j$ w.r.t. the orderly ranking by $\mathbb{B}_j$. We have $I_i$ is not reducible to zero w.r.t. $\mathbb{B}_j$ for $t+1 \Le j \Le a$. According to Lemma~\ref{LoseVar} there exists an element with the leader $u_j'$ in each $\mathbb{B}_j$ such that $u_i = \theta_{i,j}u_j'$ for some $\theta_{i,j} \in \Theta$ and all $j,$ $t+1 \Le j \Le a$.

Consider a particular $P_r$ for some $r$, $t+1 \Le r \Le a$. Let, for simplicity, $u_s = \theta_sy_s$ for all $1 \Le s \Le p$. Let $C_i$ depend on derivatives of $\tilde{Y} = \{y_{j_1}, \ldots, y_{j_{k_i}}\} \subset \{y_1,\ldots, y_{i-1}\} = Y,$ its leader $u_i = \theta_iy_i$, and some non-leaders.  Introduce a ``shifted'' function $\ord'$ as follows. For $h' = \max\limits_{1 \Le j \Le i}\ord(C_i,y_j)$, we put $$\ord'\theta y_j = \ord\theta + h' - \ord(C_i, y_j).$$
Hence, the $\ord'$-orderly ranking $>_{\ord'}$ appears on the set $\Theta\tilde{Y}$. That is, if $\ord'u > \ord'v$ then $u >_{\ord'} v$ for all $u,v \in \Theta\tilde{Y}$. Choose any rankings on the sets of derivatives $\Theta(Y\setminus\tilde{Y})$ and $\Theta \{y_{i+1},\ldots,y_l\}$. Let $\mathbb{D} = D_1,\ldots,D_{p_r}$ be a characteristic set of $P_r$ w.r.t. the ranking $$\Theta(Y\setminus\tilde{Y})>_{el}\Theta\tilde{Y} >_{el} y_i >_{el} \Theta \{y_{i+1},\ldots,y_l\}.$$ According to Theorem~\ref{OrderlyElimination} we have $\ord (D_q, y_j) \Le h$ for all $q$, $1\Le q\Le r$ and $j,$ $1 \Le j \Le l$.

The polynomial $C_i$ is reducible w.r.t. $\mathbb{D}$. Then, there exists $s,$ $1\Le s \Le p_r,$ such that one can apply a step of reduction to $C_i$ w.r.t. $D_s$. Hence, $\ord'u_{D_s} \Le \ord'$ of some variable from $C_i$. Remembering the definition of $\ord'$ and the fact that $\mathbb{C}$ is a characteristic set of $I$ w.r.t. the orderly ranking we get $\ord(D_s, y_j) \Le \ord (C_i, y_j)$ for all $j$, $1 \Le j \Le i-1$. We have {\it two} possible cases.

The first one is that $\ord(D_s, y_j) = \ord (C_i, y_j)$ for some $1 \Le j \Le i-1$. Then, either $\ord (C_i, y_i) \Le h$ and we have nothing to prove, or $\ord (C_i, y_i) = h+m> h$. Let us take care of this case. We have
$$\ord'(D_s, y_j) = \ord' (C_i, y_j) = \ord' (C_i,y_k)$$ for all $1 \Le k \Le i-1$. Let $u_{D_s} = y_{i_s}^{(h_s)}$. Since $\ord'u_{D_s} \Ge \ord'(D_s, y_j)$, we obtain that $\ord'u_{D_s} = \ord'(C_i, y_{i_s})$, that is, $h_s = \ord(C_i, y_{i_s})$.

Represent $C_i = I_{C_i}\left(y_i^{(h+m)}\right)^{n_i} + a_1\left(y_i^{(h+m)}\right)^{n_i-1} + \ldots + a_{n_i}.$ So, according to the previous paragraph one can apply an algebraic pseudo-reduction to $I_{C_i}$ and $a_q$ w.r.t. $D_s$ for all $q,$ $1 \Le q \Le n_i$.  We obtain a new reduced $C_i$ multiplied by $I_{D_s}$. Remember that $r \Ge t+1$ and, thus, both $I_{C_i}$ and $I_{D_s}$ do {\it not} lie in $P_r$ so that after this step we obtain new $\tilde{C_i}$ with the initial does not belong to $P_r$, because the ideal $P_r$ is prime. 

Since $C_i$ is reducible to zero   w.r.t. $\mathbb{D}$, we get a new $D_s$ that can reduce the polynomial $\tilde{C_i}$. If $y_i^{(h+m)}$ is a derivative of the leader of this $D_s$, we are done in this case, since the differential polynomial $D_s$ belongs to $k\{y_i,\ldots,y_l\}$ and $\ord D_s \Le h$. Taking into account that $I_{\tilde{C_i}}$ is {\it not} reducible to zero and continuing this {\it finite} process of reduction we get another $D_s$ in $P_r$ with the property $\ord(D_s, y_j) < \ord (C_i, y_j)$ for all $j$, $1 \Le j \Le i-1$. That is the second case. In this one we can put $$\widehat{C_i} = \prod\limits_{\alpha=t+1}^aD_{s_\alpha,\alpha},$$ 
where $D_{s_\alpha,\alpha}$ comes from $P_\alpha$ by the just described procedure.

Finally, multiplying $\widehat{C_i}$ by $I_{C_i}$ we obtain an element of $I$ that is not reducible w.r.t. $\mathbb{C}$. Contradiction.
Thus, we know an upper bound for the orders of $C_i$, $1\Le i \Le p$. Moreover, by Theorem~\ref{Prime} we have $\max\limits_{1 \Le j \Le a}\ord\mathbb{B}_j = \max\limits_{1\Le j \Le n}\ord\mathbb{C}_j$. By Theorem~\ref{Rosenfeld} and Lemma~\ref{LoseVar} we obtain that $$[\mathbb{C}_i]:H_{\mathbb{C}_i}^\infty \cap k[\theta_iy_i,\ord\theta_i  \Le h] = (\theta_i^jC_i^j, \ord\theta_i^j u_{C_i^j} \Le h):H_{\mathbb{C}_i}^\infty$$
and $\ord u_i \Le h$ for each $i$, $1 \Le i \Le n$. Thus, $I \cap k[\theta_iy_i,\ord\theta_i \Le h] = I'$.

Here we used the fact that the set of leaders of any characteristic set of the ideal $P$ coincides with the one of $\mathbb{C}_i$ if $P$ is a minimal prime of $[\mathbb{C}_i]:H_{\mathbb{C}_i}^\infty$. This is true due to Theorem~\ref{T20} and  Theorem~\ref{Prime}. If $P$ is a minimal prime of $I$ then $P$ is a minimal prime of some characteristic component $[\mathbb{C}_i]:H_{\mathbb{C}_i}^\infty$ of $I$ and vice versa.
Hence, our problem is purely commutative algebraic now.

In order to get a characteristic set of the ideal $I$ it is sufficient to compute an algebraic characteristic set $\mathbb{C}'$ of  $I'$ w.r.t. the {\it induced} ranking on the variables of this polynomial ring and then find
in $\mathbb{C}'$ the {\it lowest differentially autoreduced} subset $\mathbb{C}$. This set differentially reduces the ideal $I'$ to zero and some characteristic set of $I$ lies in $I'$. Thus, $\mathbb{C}$ is a characteristic set of $I$.
\end{proof}

\begin{remark} The bound obtained in Theorem~\ref{OrdinaryBound} most probably is not true for partial derivatives. At least, the method of ``shifted'' rankings does not work in non-ordinary cases. This will be shown in Example~\ref{PartialBound}.
\end{remark}

\begin{remark}\label{ComparisonSadik} One can also get another bound using \cite[Lemma 4.3]{Sadik2}. It is going to be the bound from Theorem~\ref{OrdinaryBound} multiplied by the number of differential indeterminates. For example, in $k\{y_1,\ldots,y_l\}$ it will be equal to $l\cdot h$ that is greater than $h$. So, we have got a better bound using our method.
\end{remark}

\begin{remark} By means of the proof of Theorem~\ref{OrdinaryBound} one can also complete the proof of \cite[Theorem 4.4]{Sadik2}
showing how to get an element of the ideal $I$ having elements of $P_i,$ $1 \Le i \Le a$ such that the new set becomes autoreduced. This is not clear from Sadik's paper. 
\end{remark}

\end{subsection}

\begin{subsection}{Algorithm} As an immediate consequence of Theorem~\ref{OrdinaryBound} we have the following algorithm.
 Let some orderly ranking be fixed.
\begin{algorithm}\label{CharOrd}{\sf Ordinary Characteristic Set Computation}

{\sc Input:} a finite set $F$ of ordinary differential polynomials.

{\sc Output:} characteristic set of $\{F\}$ in Kolchin's sense.

\begin{itemize}
\item Let $\mathfrak{C} =$ {\sf $\chi$-Decomposition}$(F)$ and $\mathfrak{C} = \mathbb{C}_1,\ldots,\mathbb{C}_n$ with $\mathbb{C}_i = C_i^1,\ldots,C_i^{p_i}$.
\item Let $h = \max\limits_{1\Le i \Le n}\ord\mathbb{C}_i$.
\item Compute $I' = \bigcap_{i=1}^n(\theta_i^jC_i^j, \ord\theta_i^j u_{C_i^j} \Le h):H_{\mathbb{C}_i}^\infty$.
\item $\mathbb{C}' :=$ an algebraic characteristic set of the ideal $I'$.
\item Return the differentially autoreduced subset of $\mathbb{C}'$ with the lowest rank.
\end{itemize}
\end{algorithm}

\begin{remark} The last steps of Algorithm~\ref{CharOrd} can be performed  by means of computations discussed in \cite{Wei} and \cite{Gal1, Gal2, Gal3}.
More precisely, one can compute each $I_i = (\theta_i^jC_i^j, \ord\theta_i^j u_{C_i^j} \Le h):H_{\mathbb{C}_i}^\infty$ using the Rabinovich trick and the elimination technique. Then, the intersection of the ideals $I' = \bigcap_{i=1}^n I_i$ has to be computed. The solutions to these two problems are presented in \cite{Wei}. Finally, an algorithm for computing an algebraic characteristic set of $I'$ is given in \cite{Gal1, Gal2, Gal3}.
\end{remark}

\begin{remark} Algorithm {\sf $\chi$-Decomposition} used in Algorithm~\ref{CharOrd} can be replaced by Algorithm {\sf Rosenfeld\_Groebner} presented in \cite{Bou1} and \cite{Bou2} and implemented in Maple.
\end{remark}

\end{subsection}
\end{section}

\begin{section}{Partial differential case}

\begin{subsection}{Theoretical bases}\label{Ord1}

We are {\it not} in the ordinary case now. Introduce a special class of radical differential ideals that are good in the computational sense.

\begin{definition}\label{Consistency} We say that a radical differential ideal $I$ satisfies the {\it property of consistency} iff there exists a characteristic set $\mathbb{C} \subset I$ such that $1 \notin [\mathbb{C}]:H_\mathbb{C}^\infty$.
\end{definition}

It is clear that any characterizable radical differential ideal satisfies the property of consistency. Moreover, it follows from Theorem~\ref{T20} and Theorem~\ref{Prime} that any proper regular differential ideal (an ideal of the form $[\mathbb{A}]:H_\mathbb{A}^\infty$ for a coherent autoreduced set $\mathbb{A}$) satisfies this property. Consider an example of non-regular radical differential ideal satisfying the property of
consistency.

\begin{example} Let $\mathbb{A} = x(x-1),$ $xy,$ $xz$ in $k\{x,y,z\}$ with $x < y < z$. We have $1 \notin [\mathbb{A}]:H_\mathbb{A}^\infty$. Consider
the minimal prime decomposition:
$$
\{\mathbb{A}\} = [x]\cap[x-1, y, z].
$$
We see that $\mathbb{A}$ is a characteristic set of $\{\mathbb{A}\}$. Then the radical differential ideal $\{\mathbb{A}\}$ satisfies the property of consistency. Nevertheless, since minimal primes of $\{\mathbb{A}\}$ have different sets of leader,  $\{\mathbb{A}\}$ is not a regular ideal
due to Theorem~\ref{T20} and Theorem~\ref{Prime}.
\end{example}

So, we are ready to prove Theorem~\ref{Compute}.
Let $\theta y_i$ be a differential variable in $k\{y_1,\ldots,y_l\}$. Then, by definition its order equals $\ord\theta$.

\begin{theorem}\label{Compute}  Let $I$ be a radical differential ideal in $k\{y_1,\ldots,y_l\}$ satisfying the property of consistency and a characteristic set $\mathbb{C} \subset I$ with $1 \notin [\mathbb{C}]:H_\mathbb{C}^\infty$.
\begin{enumerate}
\item
Let $U$ be the set of leaders of $\mathbb{C}$ and $U'$ be the set of leaders of any characteristic decomposition of $I$. Then $U \subset U'$.
\item Let $I = \bigcap_{i=1}^n[\mathbb{C}_i]:H_{\mathbb{C}_i}^\infty$ be a characteristic decomposition w.r.t. an orderly ranking with $\mathbb{C}_i = C_i^1,\ldots,C_i^{p_i}$. Let $h$ be the maximal order of differential variables appeared in the elements of $\mathbb{C}_i$ for $1\Le i \Le n$. Then the lowest differentially autoreduced subset of a characteristic set of the ideal $\bigcap_{i=1}^n(\theta_i^jC_i^j, \ord\theta_i^j u_{C_i^j} \Le h):H_{\mathbb{C}_i}^\infty$ is a characteristic set of $I$ w.r.t. the orderly ranking.
\end{enumerate}
\end{theorem}
\begin{proof}
We have $[\mathbb{C}] \subset I \subset [\mathbb{C}]:H_{\mathbb{C}}^\infty$. Consider the minimal prime decomposition $I = \bigcap_{i=1}^mP_i$. We have
$$
I = \bigcap_{i=1}^n[\mathbb{C}_i]:H_{\mathbb{C}_i}^\infty = \bigcap_{i=1}^mP_i.
$$

Some components of the characteristic decomposition may appear to be unnecessary.
Let $I = \bigcap_{i=1}^k[\mathbb{C}_i]:H_{\mathbb{C}_i}^\infty$ be a minimal characteristic decomposition, that is, $I \ne \bigcap_{i=1, i\ne j}^k[\mathbb{C}_i]:H_{\mathbb{C}_i}^\infty$ for all $1 \Le j \Le k$. Let $U'$ be the union of leaders of $\mathbb{C}_i$  for $1 \Le i \Le k$. If $P$ and $P_j$ are prime ideals for each $1 \Le j \Le t$ and $P \supset \bigcap_{i=1}^tP_i$ then $P \supset P_i$ for some $1 \Le j \Le t$ (see \cite[Proposition 1.11]{Atiah}).
Thus, if $P$ is a minimal prime of $I$ then $P$ is a minimal prime of $[\mathbb{C}_i]:H_{\mathbb{C}_i}^\infty$ for some $1 \Le i \Le k$.

We obtain that the set of leaders of any characteristic set of $P$ is equal to those of $\mathbb{C}_i$ by Theorem~\ref{T20} and Theorem~\ref{Prime} . Hence, the union of leaders of characteristic sets of minimal primes of $I$ is equal to $U'$. Include $[\mathbb{C}]:H_{\mathbb{C}}^\infty$ into a characteristic decomposition of $I$. For this purpose represent $\mathbb{C}$ as an output of {\sf Coherent-Autoreduced} algorithm (see \cite[Algorithm 5.1]{Fac}). Thus, we have $I = [\mathbb{C}]:H_{\mathbb{C}}^\infty \cap [\mathbb{B}_2]:H_{\mathbb{B}_2}^\infty \cap \ldots \cap [\mathbb{B}_r]:H_{\mathbb{B}_r}^\infty$. Denote the set of leaders of $\mathbb{C}$ by $U$.

Let $P$ be a minimal prime of $[\mathbb{C}]:H_{\mathbb{C}}^\infty$. Then $P$ is a minimal prime of $\{\mathbb{C}\}$ (see \cite[page 644]{Fac}). Thus, $P$
is a minimal prime of $I$, because $\{\mathbb{C}\} \subset I \subset [\mathbb{C}]:H_{\mathbb{C}}^\infty$. Since $P$ is a minimal prime of $[\mathbb{C}_i]:H_{\mathbb{C}_i}^\infty$ for some $1 \Le i \Le k$, then the set of leaders of $\mathbb{C}_i$ is equal to $U$ and $U \subset U'$. So, we know an upper bound for the order of a characteristic set of $I$.  Due to Theorem~\ref{Rosenfeld} and Lemma~\ref{LoseVar} we have $[\mathbb{C}_i]:H_{\mathbb{C}_i}^\infty \cap k[\theta_iy_i,\ord\theta_i  \Le h] = (\theta_i^jC_i^j, \ord\theta_i^j u_{C_i^j} \Le h):H_{\mathbb{C}_i}^\infty$ and $\ord u_i \Le h$ for each $1 \Le i \Le n$. Thus, we obtained the result because the end of the proof is the same as in Theorem~\ref{OrdinaryBound}.
\end{proof}

The main contribution of Theorem~\ref{Compute} is that our computations are moved into the ring of commutative polynomials in a {\it finite} number of variables. This is a crucial point in Algorithm~\ref{CharAlg}.

\begin{remark}\label{Concentrated} We see that the set of leaders $U$ of $\mathbb{C}$ is not only a subset of $U'$. We have $U$ is {\it equal} to the set of leaders of some characteristic component in any characteristic decomposition of $I$. Thus, $U$ is ``concentrated'' in some characteristic component. We call it the ``localization'' property.
\end{remark}

\begin{remark}\label{ConcentratedArbitrary} The fact that the set of leaders of $\mathbb{C}$ is equal to that of some characteristic component holds true for {\it any} differential  ranking.
This follows from the proof of Theorem~\ref{Compute}.
\end{remark}

Look at the following example of a {\it non-regular} radical differential ideal
illustrating Theorem~\ref{Compute}.

\begin{example} Consider the radical differential ideal $I = \{xy\}$ in $k\{x,y\}$ w.r.t.
any differential ranking such that $x < y$. We have the following decomposition:
$$
I = [x]\cap[y].
$$
The ideal $I$ is not prime. Moreover, it is not neither characterizable nor regular.
Nevertheless, the ideal $I$ satisfies the property of consistency, because one
can consider its characteristic set $\mathbb{C} = xy$ in Kolchin's sense with $1 \notin [xy]:x^\infty$.

As Theorem~\ref{Compute} tells us the set of leaders $y$ of the characteristic set $\mathbb{C}$ must equal to the set of leaders of some component of a characteristic decomposition of the ideal $I$. This is the case for its second component $[y]$.
In order to compute $\mathbb{C}$ we do not need to differentiate anything. We just
calculate the reduced Gr\"obner basis of $(x)\cap(y)$, that is equal to $xy$, and
output this set as a characteristic set of the ideal $I$, since $xy$ is a differentially
autoreduced set.
\end{example}

Note that the condition that an ideal satisfies the property of consistency cannot be omitted in Theorem~\ref{Compute}. To support this fact we give Example~\ref{CounterExample}.

\end{subsection}

\begin{subsection}{Algorithm}\label{Ord2}

In conclusion, we obtain the following algorithm. Let some orderly ranking be fixed.
\begin{algorithm}\label{CharAlg}{\sf Characteristic Set Computation}

{\sc Input:} a finite set $F$ of differential polynomials such that $\{F\}$ satisfies the property of consistency.

{\sc Output:} characteristic set of $\{F\}$ in Kolchin's sense.

\begin{itemize}
\item Let $\mathfrak{C} =$ {\sf $\chi$-Decomposition}$(F)$ and $\mathfrak{C} = \mathbb{C}_1,\ldots,\mathbb{C}_n$ with $\mathbb{C}_i = C_i^1,\ldots,C_i^{p_i}$.
\item Let $h = \max\limits_{1\Le i \Le n}\max\limits_{1 \Le j \Le p_i} \ord C_i^j$.
\item Compute $I' = \bigcap_{i=1}^n(\theta_i^jC_i^j, \ord\theta_i^j u_{C_i^j} \Le h):H_{\mathbb{C}_i}^\infty$.
\item $\mathbb{C}' :=$ an algebraic characteristic set of the ideal $I'$.
\item Return the differentially autoreduced subset of $\mathbb{C}'$ with the lowest rank.
\end{itemize}
\end{algorithm}

The last steps of Algorithm~\ref{CharAlg} can be performed  in the same way as done in Algorithm~\ref{CharOrd}.

\begin{remark} Algorithm {\sf $\chi$-Decomposition} used in Algorithm~\ref{CharAlg} can be replaced by Algorithm {\sf Rosenfeld\_Groebner} presented in \cite{Bou1} and \cite{Bou2} and implemented in Maple.
\end{remark}

Note that $h$ from Algorithm~\ref{CharAlg} is less than $h$ from Algorithm~\ref{CharOrd} in general. Hence, the ideals satisfying the property of consistency are {\it better} in the computational sense than radical differential ideals in general.

\end{subsection}

\begin{subsection}{Arbitrary differential rankings}\label{Arbitrary}
Theorem~\ref{Compute} can be easily applied to computation of Kolchin characteristic sets of radical differential ideals satisfying the property of consistency {\it not only} in the orderly case. The statement of the theorem remains the same  eliminating the word ``orderly''. But in order to move the problem to a commutative ring in a finite number of variables we also need to use a characteristic decomposition w.r.t. an orderly ranking. 

So, the algorithms coming from the proof of Theorem~\ref{Compute} is the following. Let some differential ranking $>$ be fixed.
\begin{algorithm}\label{CharAlgArbitrary}{\sf Characteristic Set Computation}

{\sc Input:} a finite set $F$ of differential polynomials such that $\{F\}$ satisfies the property of consistency.

{\sc Output:} characteristic set of $\{F\}$ in Kolchin's sense.

\begin{itemize}
\item Let $\mathfrak{B} =$ {\sf $\chi$-Decomposition}$(F)$ w.r.t. the ranking $>$ and $\mathfrak{B} = \mathbb{B}_1,\ldots,\mathbb{B}_k$ with $\mathbb{B}_i = B_i^1,\ldots,B_i^{p_i}$.
\item Let $h = \max\limits_{1\Le i \Le k}\max\limits_{1 \Le j \Le p_i} \ord B_i^j$.
\item Let $\mathfrak{C} =$ {\sf $\chi$-Decomposition}$(F)$ w.r.t. an orderly ranking and $\mathfrak{C} = \mathbb{C}_1,\ldots,\mathbb{C}_n$ with $\mathbb{C}_i = C_i^1,\ldots,C_i^{q_i}$.
\item Compute $I' = \bigcap_{i=1}^n(\theta_i^jC_i^j, \ord\theta_i^j u_{C_i^j} \Le h):H_{\mathbb{C}_i}^\infty$.
\item $\mathbb{C}' :=$ an algebraic characteristic set of $I'$ w.r.t. the ranking induced by $>.$
\item Return the differentially autoreduced subset of $\mathbb{C}'$ with the lowest rank.
\end{itemize}
\end{algorithm}

\begin{corollary} There exists an algorithm for computing a characteristic set of a radical differential ideal satisfying the property of consistency in the partial differential case.
\end{corollary}

\end{subsection}

\end{section}

\begin{section}{Examples}
We show how to apply Algorithm~\ref{CharAlg} to a particular radical differential ideal.

\begin{example} Let $I = \{(x-t)x',$ $x'y',$ $(x-t)(z'+y')\}$ in $\mathbb{Q}(t)\{x,y,z\}$ with $t' = 1$ and an orderly ranking $x < y < z$. We have the following decomposition:
$$
 I = [x-t, y']\cap[x',z'+y'].
$$

So, the maximal order of variables appeared in this decomposition is equal to $1$. Hence, we need to compute the reduced Gr\"obner basis $G$ of the ideal $I' = (x-t, x'-1,y')\cap(x',z'+y')$. This can be done by the elimination technique: $G$ equals the intersection of the reduced Gr\"obner basis  w.r.t. the lexicographic ordering $x < x' < y' < z' < w$ of the ideal $(w(x-t),$ $w(x'-1),$ $wy',$ $(1-w)x',$ $(1-w)(z'+y'))$ and the ring $\mathbb{Q}(t)[x, x', y', z']$.

Finally, $G = (x-t)x',$ $x'(x'-1),$ $x'y',$ $(x-t)(z'+y'),$ $x'z'-(z'+y'),$ $y'(z'+y')$. Then a characteristic set of $I'$ equals $\mathbb{C} = (x-t)x',$ $(x-t)(z'+y')$ and by Theorem~\ref{Compute} a characteristic set of the radical differential ideal $I$ is also $\mathbb{C}$.
\end{example}

The following example shows the {\it difference} between radical differential ideals with the consistency property and an arbitrary radical
differential ideal. This can be considered as a ``counter-example'' for Remark~\ref{Concentrated} and Remark~\ref{ConcentratedArbitrary}.
\begin{example}\label{Localization} Let $\mathbb{A} = x(x-1),$ $xy,$ $(x-1)z$. We have $1 \in [\mathbb{A}]:H_\mathbb{A}^\infty$. Consider
the minimal prime decomposition:
$$
\{\mathbb{A}\} = [x, z]\cap[x-1, y].
$$
We see $\mathbb{A}$ is a characteristic set of $\{\mathbb{A}\}$ with the set of leaders equals $U = x,$ $y,$ $z$.
Thus, $U$ does not necessarily correspond to the {\it unique} characteristic component. In this example $U \subset x,z \cup x,y$ and the localization property is not valid.
\end{example}

Note that Theorem~\ref{Compute} is true for Example~\ref{Localization}: in this case we have a restriction to the orders of the elements of a characteristic set $\mathbb{A}$.
Consider a ``counter-example'' for Theorem~\ref{Compute}.
\begin{example}\label{CounterExample} Consider a radical differential ideal defined by its characteristic decomposition:
$$
I = [x-1, y]\cap[x, y^{(n)}, z^{(m)}+y]
$$
in $k\{x, y, z\}$ with $x < y < z$, an orderly ranking, and $n \Le m$. Both of these components are prime differential ideals, because they are generated by linear differential polynomials. In addition, since they are prime, these radical differential ideals are also characterizable (see \cite[page 646]{Fac}). One can show that a characteristic set of $I$ is $\mathbb{C} = x(x-1),$ $xy,$ $(x-1)z^{(m+n)}.$ The radical differential ideal $I$ does not satisfy the property of consistency and Theorem~\ref{Compute} is not true
for $I$. Indeed, for $m, n > 0$ we have $m+n > \max\{m, n\}.$
\end{example}

So, we see that the upper bound established in Theorem~\ref{Compute} is {\it wrong} for some radical differential ideals {\it not satisfying} the property of consistency. Nevertheless, Example~\ref{CounterExample} is in the ordinary case and Theorem~\ref{OrdinaryBound} can be applied. We have $\ord \{x, y^{(n)}, z^{(m)}+y\} = m + n$ and the maximal order of elements of $\mathbb{C}$ does not exceed $m+n$.

\begin{remark} We see that the bound for radical differential ideals satisfying the property of consistency is lower than one for arbitrary radical differential ideals in general. That is why radical differential ideals satisfying the property of consistency are better in the computational sense. Indeed, Algorithm~\ref{CharAlg} works faster than Algorithm~\ref{CharOrd}, since we are in a lower algebraic dimension in Algorithm~\ref{CharAlg} than in Algorithm~\ref{CharOrd}.
\end{remark}

\end{section}

\begin{section}{Conjecture}\label{Conjecture}
The following example shows that the method of ``shifted'' rankings used in Theorem~\ref{OrdinaryBound} does not work in the partial derivative case.

\begin{example}\label{PartialBound}
Consider $k\{u\}$ with $\Delta = \{\partial/\partial x, \partial/\partial y,\partial/\partial z\}$.
Let $$P = [u_{yyyy}+u_{xxyy}+u_{xxxx}+u_z,\ u_{xxxxx}+u_{xxy}+u_{xyy}-u_{xxyyy}].$$
Consider $f = u_{xxy}+u_{xyy}+u_{xxxxx}+u_{yyyyy}+u_{xxxxy}+u_{yz} \in P.$ One can prove that this polynomial cannot be an element of any characteristic set of $P$ w.r.t. ``shifted'' rankings. However, $f$ can be an element of a characteristic set of a radical differential ideal for which the ideal $P$ is a minimal prime component.
\end{example}

\begin{remark} Example~\ref{PartialBound} is not a counter-example for Theorem~\ref{OrdinaryBound} in non-ordinary cases. It is just a counter-example for the part of the proof of the theorem.
\end{remark}

\begin{corollary}
One need to develop another technique to get a bound for the order of a characteristic set of a radical differential ideal in partial derivative cases and radical differential ideals not satisfying the property of consistency.
\end{corollary}

Nevertheless, we hope that this problem can be solved by means of increasing estimates of bounds for characteristic sets of radical differential ideals. So, we {\it conjecture} that if $|\Delta| = m$ then the order of each element of a characteristic set of a radical differential ideal $I$ is bounded by $qh^m$, where $h = \max\limits_{1\Le i \Le n}\ord\mathbb{C}_i$ and $I = \bigcap_{i=1}^n[\mathbb{C}_i]:H_{\mathbb{C}_i}^\infty$ is a characteristic decomposition w.r.t. an {\it orderly} ranking and $q$ is the number of different differential indeterminates $y_j$ appearing as leaders in $\mathbb{C}_i$ for $1 \Le i \Le n$.
\end{section}

\begin{section}{Conclusions}
We discussed an algorithm for computing a Kolchin characteristic set of an arbitrary ordinary radical differential ideal w.r.t. orderly rankings.
In the partial differential case we also presented a solution to the problem of computing a characteristic set of a radical differential ideal satisfying the special property
of consistency. These solutions are new and previously these problems were completely solved only in the non-differential case as to our knowledge.

The authors hope that the technique obtained in this paper can be generalized to any radical differential ideals using the ideas we presented. Another natural way of generalizing
these results is to complete our investigation of non-orderly rankings such as, for example, very important elimination ones in the case of partial differential polynomials.
\end{section}

\begin{section}{Acknowledgements}
We thank Evgeniy V. Pankratiev and Michael F. Singer for their helpful comments and support. We are appreciate to Alexey Zobnin, Alexander B. Levin, and Oleg Golubitsky
for their reviewing of the paper and for valuable suggestions. We also thank referees
of the ACA 2004 for their important remarks.
\end{section}

\end{document}